\documentstyle[12pt]{article}
\newcommand{\const}{\mathop{\rm const}\limits}

\newcommand{\supp}{\mathop{\rm supp}\limits}

\newcommand{\Var}{\mathop{\rm Var}\limits}

\begin{document}

\begin{center}

{\bf  LOWER BOUNDARIES FOR PARAMETRIC ESTIMATIONS \\

\vspace{3mm}

IN DIFFERENT NORMS} \par

\vspace{4mm}

 $ {\bf E.Ostrovsky^a, \ \ L.Sirota^b } $ \\

\vspace{3mm}

$ ^a $ Corresponding Author. Department of Mathematics and computer science, Bar-Ilan University, 84105, Ramat Gan, Israel.
E-mail: \  eugostrovsky@list.ru\\

\vspace{3mm}

$ ^b $  Department of Mathematics and computer science. Bar-Ilan University,
84105, Ramat Gan, Israel. \ E-mail: sirota3@bezeqint.net \\

\vspace{4mm}

{\bf Abstract.} \\

 \end{center}

\vspace{3mm}

 We establish some new {\it non-asymptotical} lower bounds for deviation of regular unbiased estimation of unknown
 parameter from its true value in different norms, alike the classical Rao - Kramer's  inequality.\par
 We show that if the new norm is weaker that ordinary Hilbertian  norm, that the rate of convergence of arbitrary regular unbiased estimate
 does not exceed $ 1/\sqrt{n}, $  and if the new norm is stronger that one, the rate of convergence of the well-known Maximal Likelihood
 Estimate (MLE) is also equal to $ 1/\sqrt{n}. $\par

\vspace{3mm}

 \par

\vspace{4mm}

\section{  Statement of problem. Notations. Assumptions. }

\vspace{3mm}

 Let $ (\Omega, B, {\bf P} ) $ be non - trivial probability space with expectation $  {\bf E } $  and variance  $ \Var, $
 $  (X = \{x \}, \cal{A}, \mu) $ be measurable space equipped
with sigma - finite measure $  \mu, \ \Theta = \{  \theta \} $    be connected subset of real line, i.e. open, semi-open or closed interval,
$ \theta_0 $ be a fixed interior point in the set $ \Theta.  $ \par
 It is sufficient to suppose for example building that on the set $ \Omega $ there exists an uniform distributed random variable.\par

 Let also $ f = f(x,\theta), \ x \in X, \ \theta \in \Theta $ be differentiable relative the parameter $ \theta $ strictly positive density, i.e.
numerical measurable normed function:

$$
\int_X f(x,\theta) \ \mu(dx) = 1, \ \forall x \in X, \ \forall \theta \in \Theta \ \Rightarrow f(x,\theta) > 0..\eqno(1.1)
$$

 We suppose that the random variable $ \xi: \Omega \to R  $ has a density $  f(x,\theta_0): $

$$
{\bf P}(\xi \in G) = \int_G f(x,\theta_0) \ \mu(dx),
$$
i.e. the value $ \theta_0 $ is true value of the parameter $ \theta. $\par

 Let us denote by $  L(\xi) = L(\xi, \theta) $ the ordinary likelihood function:

$$
L(\xi) = L(\xi,\theta) \stackrel{def}{=} \log f(\xi,\theta).
$$
 The following function

$$
\theta \to L(\xi,\theta) - L(\xi, \theta_0) = \log [ f(\xi, \theta)/f(\xi,\theta_0)  ]
$$
 is named  a {\it  contrast} function. \par

 The r.v. $  \xi  $ may be also a vector with values  in  the space $ R^n, $ in particular, may be a sample of a volume $  n: $

$$
\xi = \vec{\xi} = \{ \eta(1), \eta(2), \ldots, \eta(n)   \}, \ n = 2,3,\ldots;
$$
 where  the r.v.  $ \{ \eta(i)  \}  $ are i. i.d. with the positive density which we denote by $ g(x, \theta). $ Of course,

$$
f(x,\theta) = \prod_{i=1}^n g(x_i, \theta), \ x = \{ x(1), x(2), \ldots, x(n)  \}.
$$
 We denote also in this case

$$
l_i = l_i(\eta(i)) = l_i(\eta(i), \theta) = \partial \log g(\eta(i), \theta)/\partial \theta; \ l = l_1.
$$

  Further, let $ \hat{\theta} = \hat{\theta}(\xi) $ be some unbiased {\it regular} in the sense of the monograph
\cite{Ibragimov1}  estimate of the parameter $ \theta: $

$$
{\bf E} \hat{\theta} = \int_X \hat{\theta}(x) \ f(x,\theta) \ \mu(dx) = \theta, \ \theta \in \Theta, \eqno(1.2)
$$
 All we need is-the following two equalities:

$$
{\bf E} \left[(\hat{\theta} - \theta_0) \cdot \frac{\partial \log f(\xi,\theta)}{\partial \theta} \right] =
\int_X \left[ (\hat{\theta}(x) - \theta_0) \cdot \partial \log f(x,\theta)/\partial \theta \right] \ \mu(dx)=1; \eqno(1.3)
$$

$$
{\bf E}  \frac{\partial \log f(\xi,\theta)}{\partial \theta} =
\int_X  \left[ \partial \log f(x,\theta)/\partial \theta \right] \ \mu(dx) =0. \eqno(1.4)
$$

\vspace{4mm}

{\bf  Our purpose in this report is obtaining  the lower non - asymptotical estimation of Rao-Kramer's type for the deviation
$  \sqrt{n} ||(\hat{\theta}_n - \theta_0)||Z, $ i.e. under the classical norming sequence $ \sqrt{n}, $
 for some different r.i. norms over source probability space $ ||\cdot||Z. $  }\par

\vspace{3mm}

The {\it  upper } non - asymptotical estimations for these deviation, and as a consequence an exponentially  exact confidential interval
for the unknown parameter $ \theta_0, $
under modern terms: majorizing measures, generic chaining etc. for the MLE estimates was derived in the article
\cite{Ostrovsky1}; see also \cite{Borovkov1},  chapter 2, section 23; \cite{Korostelev1},  chapter 3, Lemma 3.19. \par

\vspace{3mm}

 It is clear that the norm $ ||\cdot||Y $ should be substantially weaker as the classical $ L_2(\Omega, {\bf P}), $ as well as
for the investigation of upper estimation this norm should be stronger as one. \par

\vspace{3mm}

 Note briefly that the case of biased estimate,  multivariate parameter and both this circumnutations may be investigated
quite analogously.\par

\vspace{4mm}

\section{ General estimates }

\vspace{3mm}

 Let $  (Y, ||\cdot||Y) $ be arbitrary rearrangement invariant (r.i.) space over $ (\Omega, B, {\bf P} ). $  Reader can found using for us
 facts about the theory of this spaces in the classical monograph \cite{Bennet1}.\par

 {\it We accept that all considered in this article r.i. spaces will be constructed over our probability  space } $ ( \Omega, B, {\bf P}). $ \par

 We denote  as ordinary by $ (Y', ||\cdot||Y') $ the associate space relative the "scalar product"

 $$
 ( \zeta, \tau) = {\bf E} \zeta \tau = \int_{\Omega} \zeta(\omega) \ \tau(\omega) \ {\bf P} (d \omega), \eqno(2.1)
 $$
so that $ (Y', ||\cdot||Y') $ is  again r.i. space and

$$
||\tau||Y' = \sup_{ \zeta: ||\zeta||Y = 1  } \left[ \frac{|(\zeta, \tau)|}{||\zeta||Y} \right]; \hspace{5mm}
||\zeta||Y = \sup_{ \tau: ||\tau||Y' = 1  } \left[ \frac{|(\zeta, \tau)|}{||\tau||Y'} \right]. \eqno(2.2)
$$
 Therefore,

$$
|(\zeta, \tau)| \le ||\zeta||Y \cdot ||\tau||Y', \eqno(2.3)
$$
the generalized H\"older's inequality.\par

 Let $ \eta $  be any centered r.v. belonging to the r.i. space $ Y. $ We define  (and denote)
 analogously  M.Ledoux and M.Talagrand \cite{Ledoux1}, p.274 - 275
 the following so-called $ CLT(Y) = C(Y) $  norm for $  \eta $ as follows:

$$
||\eta||CLT(Y) = ||\eta||C(Y) \stackrel{def}{=} \sup_n || n^{-1/2} \sum_{i=1}^n \eta_i  ||Y, \eqno(2.4)
$$
where $ \{  \eta_i \} $ are independent copies of $  \eta. $ It will be presumed  without loss of generality
that the probability space is sufficiently rich, see beginning of this article.\par
 Obviously, if $ ||\eta||C(Y) < \infty $  then $  ||\eta||L_2(\Omega, P) < \infty.  $\par

 Note but that our definition does not coincides with the definition  of M.Ledoux and M.Talagrand.\par

\vspace{3mm}

{\bf  Theorem 2.1. }  Let $ \xi = \{  \eta(i) \}, \ i = 1,,2,\ldots,n $ be a sample of a volume $ n. $
 Suppose that there exists a r.i. space $ (Y, ||\cdot||) $ over  our probability space such that the
r.v. $ l = l(\eta_1, \theta) $ belongs to the space $ CLT(Y): \ 0 < ||l||CLT(Y) < \infty. $ Then for arbitrary
unbiased regular estimate $ \hat{\theta} $

$$
\sqrt{n} ||\ \hat{\theta} - \theta_0 \ ||Y' \ge \frac{1}{||l||CLT(Y)}. \eqno(2.5)
$$

\vspace{3mm}

{\bf Proof.} We start from the relations (1.3) and (1.4):

$$
1 = (\hat{\theta} - \theta_0, \ \sum_{i=1}^n l_i). \eqno(2.6)
$$
 We use the (generalized) H\"older's inequality (2.3):

$$
1 \ge || \ \hat{\theta} - \theta_0 \ ||Y'  \cdot || \sum_{i=1}^n l_i ||Y =  \sqrt{n} \cdot
|| \  \hat{\theta} - \theta_0 \ ||Y' \cdot || n^{-1//2} \ \sum_{i=1}^n l_i ||Y.
$$

 It follows from the direct definition of the $ CLT(Y) $ norm

$$
1 \ge  \sqrt{n} \cdot || \  \hat{\theta} - \theta_0 \ ||Y' \cdot ||  l ||CLT(Y),
$$
which is equivalent to the assertion (2.5) of theorem 2.1. \par

 If for example $  Y = Y' = L_2(\Omega, {\bf P}), $ we get to the classical  inequality of Rao - Kramer. \par

\vspace{3mm}

{\bf Remark 2.1.} Note that in general case the quantity $ || \ l \ ||CLT(Y)  $ dependent on the parameter $ \theta. $ \par

\vspace{3mm}

{\bf Definition 2.1.}  The r.i. space $ (Y, ||\cdot||) $ is said to be {\it strong normal rearrangement invariant,} briefly, s.n.r.i.,
write $  Y \in s.n.r.i., $ if the auxiliary space $  (CLT(Y), ||\cdot||) $ is  equivalent to source space $ (Y, ||\cdot||) $  on the
subspace of the centered variables $ \{  \eta \} $ from this space:

$$
\forall \{ \eta(i), \ \eta(i) \in Y, \ {\bf E} \eta(i) = 0 \} \ \Rightarrow ||\sum_{i=1}^n \eta(j)||Y \le K(Y)
\sqrt{ \sum_{i=1}^n (||\eta(i)|| Y)^2 }  \eqno(2.7)
$$
for some finite constant   $  K(Y) $ depending only on the whole space $  Y. $ \par
 It is clear that for mean zero variable  $  \eta \ ||\eta||Y \le ||\eta||CLT(Y),  $ therefore in s.n.r.i. spaces both the  norm are really (linear)
equivalent:

$$
||\eta||Y \le ||\eta||CLT(Y) \le K(Y) ||\eta||Y.
$$

 The symbol  $  "n" $ in the abbreviate of definition 2.1  comes from the word "normal". \par

 \vspace{3mm}

{\bf Proposition 2.1.} If in addition to the conditions of theorem 2.1 the space $ Y $  is s.n.r.i., then evidently

$$
\sqrt{n} ||\ \hat{\theta} - \theta_0 \ ||Y' \ge \frac{1}{K(Y) \ ||l||Y}. \eqno(2.8)
$$

\vspace{3mm}

{\bf Definition 2.2.}  The r.i. space $ (Y, ||\cdot||) $ is said to be {\it weak normal rearrangement invariant,} briefly, w.n.r.i.,
write $  Y \in w.n.r.i., $ if for arbitrary centered r.v. $  \eta $  from this space the norm $ ||\eta||CLT(Y)  $ is finite.\par

\vspace{3mm}

{\bf Proposition 2.2.} If in addition to the conditions of theorem 2.1 the space $ Y $  is w.n.r.i., then

$$
\overline{\lim}_{n \to \infty}   \sqrt{n} ||\ \hat{\theta} - \theta_0 \ ||Y' > 0. \eqno(2.9)
$$

\vspace{3mm}

 Recall, see e.g. \cite{Ibragimov2}, chapters  2-4, that the centered (or moreover symmetrically distributed)  r.v. $  \eta $ belongs to the
Domain of Stable Attraction $ (DSA), $ iff

$$
x \to \infty \Rightarrow {\bf P} (\eta < -x) \sim C_1 x^{-\alpha} L(x), \ {\bf P} (\eta > x) \sim C_2 x^{-\alpha} L(x), \eqno(2.10)
$$
where $ C_1, C_2 = \const > 0, \ L(x) $ is continuous non - negative slowly varying as $ x \to \infty $ function, $ \alpha = \const \in (0,2). $\par

 We define also the set $  DNA_{\infty} $ as a set of all centered random variables $  \{ \zeta \} $  from the Domain of Normal  Attraction $ DNA $
but with finite variation:  $ \Var(\zeta) = \infty. $ \par

\vspace{3mm}

{\bf Proposition 2.3.} If the r.i. space $ (Y, ||\cdot||) $ is such that

$$
Y \cap DSA \ne \emptyset \eqno(2.11)
$$
or

$$
Y \cap DNA_{\infty} \ne \emptyset,  \eqno(2.12)
$$
then this space $ (Y, ||\cdot||) $ is not weak normal rearrangement invariant space.\par

\vspace{3mm}

  {\bf Proof.}  It is sufficient to take in the definition (2.2) the mean zero random variable $  \eta $ from the set

$$
\eta \in (Y \cap DSA) \cup ( Y \cap DNA_{\infty})
$$
to ensure that the definition (2.2) is not fulfilled.\par

\vspace{4mm}

\section{ Lebesgue - Riesz norm }

\vspace{3mm}

 We consider in this section the case when at the capacity of the space $  Y  $ is the classical LebesgRiesz space
 $ L_p = L_p(\Omega, {\bf P}), \ p = \const \in [1, \infty).  $ We  will denote as usually

$$
|\eta|_p = \left[ {\bf E} |\eta|^p  \right]^{1/p};  \hspace{5mm} q = p/(p-1), p > 1; \hspace{5mm} q = \infty, p = 1.
$$
 Define also the so - called $  p \ -$ Fisher's information $ i_p(\theta): $

$$
i_p(\theta) = i_p(\eta, \theta):= |  l(\eta, \theta) |_p = | \partial \log  g(\eta, \theta) /\partial \theta |_p, \eqno(3.1)
$$
and analogously for the sample

$$
I_p(\theta) = I(\vec{\xi}, \theta) := | l(\vec{\xi}, \theta) |_p. \eqno(3.1a)
$$
 which coincides with the classical Fisher's information when $  q = p = 2. $ \par

 The expression for $ i_p(\vec{\xi}, \theta) $ may be rewritten as follows

$$
i_p(\eta, \theta) = \int_X |g'_{\theta}(x,\theta)|^p \cdot g^{1 -p}(x,\theta) \ \mu(dx), \eqno(3.1b)
$$
and analogously for the quantity $ I_p(\vec{\xi}, \theta). $ \par

 \vspace{3mm}

{\bf Theorem 3.1.} \par

{\bf A. } The space  $ L_p = L_p(\Omega, {\bf P}) $ is not weak normal rearrangement invariant if $ 1 \le p < 2. $\par

\vspace{3mm}

{\bf B.} The space  $ L_p = L_p(\Omega, {\bf P}) $ is  strong normal rearrangement invariant if $  p \ge 2. $\par

\vspace{3mm}

{\bf Proof.}  The first proposition follows immediately from theorem 2.1, as long as the space $ L_p, \ 1 \le p < 2 $ contains the
symmetric stable distributed random variable $ \zeta  $ with the parameter $ \alpha = (p + 2)/2; \ \alpha \in (0.2): $

$$
{\bf E} e^{i t \zeta} = e^{ - |t|^{\alpha} }, \ t \in R.
$$

 So, let now $  p \ge 2.  $  The using for us inequality (2.7) is a particular case of the famous Rosenthal's inequality \cite{Rosenthal1}:

$$
|\sum_{j=1}^n l_j|_p \le R(p) \cdot \sqrt{n} \cdot |l|_p, \ p \ge 2, \eqno(3.2)
$$
where the Rosenthal's "constant" $  R(p) $ may be estimated as follows:

$$
R(p) \le C \cdot \frac{p}{e \ln p}, \hspace{5mm} C \le 1.77638\ldots \eqno(3.3)
$$
see \cite{Ostrovsky2}.\par

 Therefore, we can accept in (2.7) - (2.8) $ K(L_p) = R(p) $ with estimation (3.3).\par

 We  conclude for the regular sample and  regular unbiased estimate $ \hat{\theta} = \hat{\theta}_n: $ \\

{\bf  Proposition 3.1. }

$$
\sqrt{n} \ ||\hat{\theta}_n - \theta_0||_q \ge  \frac{1}{R(p) \ i_p(\theta_0)}, \ q \in (1,2), \ p = q/(q-1). \eqno(3.4)
$$

\vspace{3mm}

{\bf Remark 3.1.} It follows from the triangle inequality that if $  \{ \eta_i \} $ are independent, then

$$
I_p(\theta) \le \sum_{i=1}^n i_p(\eta_i, \theta),
$$
but it follows from the Rosenthal's inequality more exact as $ n >> 1 $ estimate

$$
I_p(\theta) \le R(p) \cdot \sqrt{ \sum_{i=1}^n i^2_p(\eta_i, \theta)}. \eqno(3.5)
$$

\vspace{3mm}

{\bf Remark 3.2.} The notion of  $  p \ -$ Fisher's information $ i_p(\theta) $ in (3.1) and (3.1a) may be generalized
on  arbitrary r.i. space $ (Y, ||\cdot||Y): $

$$
i_{(Y)}(\theta) = i_{(Y)}(\eta, \theta) \stackrel{def}{=}
||  l(\eta, \theta) ||Y = || \partial \log  g(\eta, \theta) /\partial \theta  ||Y, \eqno(3.6)
$$
and analogously for the sample

$$
I_{(Y)}(\theta) = I_{(Y)}(\vec{\xi},\theta)   := || l(\vec{\xi}, \theta) ||Y. \eqno(3.6a)
$$

 Obviously, if the space $ (Y, ||\cdot||Y) $ is w.n.r.i., then for independent variables (observations) $ \{  \eta(i) \} $

$$
I_{CLT(Y)}(\vec{\xi}, \theta) \le K(Y) \sqrt{ \sum_{i=1}^n  \left[ i_{(Y)}(\eta(i), \theta) \right]^2 }. \eqno(3.7)
$$

 If in addition  the space $ (Y, ||\cdot||Y) $ is s.n.r.i., then for independent variables (observations) $ \{  \eta(i) \} $

$$
I_{(Y)}(\vec{\xi}, \theta) \le K(Y) \sqrt{ \sum_{i=1}^n  \left[ i_{(Y)}(\eta(i), \theta) \right]^2 }. \eqno(3.8)
$$

\vspace{4mm}

\section{ Grand Lebesgue space norm }

\vspace{3mm}

 Recently, see \cite{Fiorenza1}, \cite{Fiorenza2},\cite{Ivaniec1}, \cite{Ivaniec2}, \cite{Jawerth1},
\cite{Karadzov1}, \cite{Kozachenko1}, \cite{Liflyand1}, \cite{Ostrovsky3} etc.
 appear the so-called Grand Lebesgue Spaces (GLS)

 $$
 G(\psi) = G = G(\psi ; (1,B));  \ B = \const \in (1, \infty]
 $$
spaces consisting on all the random variables  (measurable functions) $ f : \Omega \to R  $ with finite norms

$$
||f||G(\psi) = ||f||G(\psi; (1,B))  \stackrel{def}{=} \sup_{p \in (A;B)} \left[\frac{|f|_p}{\psi(p)} \right]. \eqno(4.1)
$$

 Here $ \psi = \psi(p), \ p \in [1,B) $ is some continuous positive on the {\it open} interval $ (1;B) $ function such
that

$$
\inf_{p \in(A;B)} \psi(p) > 0. \eqno(4.2)
$$

We will denote
$$
\supp(\psi) \stackrel{def}{=} [1;B)
$$
 or by abuse of laanguage $ \supp(\psi) = B.  $ \par

The set of all such a functions with the support $ \supp(\psi) = (1;B) $ will be denoted by  $  \Psi(1;B) = \Psi(B). $  \par

 This spaces are rearrangement invariant; and are used, for example, in the theory of Probability, theory of Partial Differential Equations,
 Functional Analysis, theory of Fourier series, Martingales, Mathematical Statistics, theory of Approximation  etc. \par

 Notice that the classical Lebesgue - Riesz spaces $ L_p $  are extremal case of Grand Lebesgue Spaces, see
 \cite{Ostrovsky3}. \par

 Let a function $  \xi:  \Omega \to R  $ be such that

 $$
 \exists B > 1  \Rightarrow  \forall p \in [1, B)  \ |\xi|_p < \infty.
 $$
Then the function $  \psi = \psi_{\xi}(p) $ may be {\it naturally} defined by the following way:

$$
\psi_{\xi}(p) := |\xi|_p, \ p \in [1,B). \eqno(4.3)
$$

  More generally, let $ \xi(\alpha), \ \alpha \in A, \  A $  is arbitrary set, be a {\it family } of mean zero r.v. such that

$$
\exists B > 2, \ \forall p \in [2,B) \Rightarrow \psi^{(A)}(p) := \sup_{ \alpha \in A} |\xi(\alpha)|_p < \infty.
$$
The function $ p \to \psi^{(A)}(p) $ is called a {\it natural } function for the family  $  \xi(\alpha), \ \alpha \in A. $ \par
We emphasize that the variables $ \{  \xi(\alpha) \} $ ​​can be arbitrarily dependent and that

$$
\sup_{\alpha \in A} ||\xi(\alpha)||G(\psi^{(A)}) = 1.
$$

 The finiteness of the $ G\psi \ - $ norm for some r.v. $  \xi $ allows to obtain the exact exponential tail inequalities for the
distribution  $  \xi; $ for instance,

$$
\sup_{p \ge 1} \left[ \frac{|\xi|_p}{p^{1/m}} \right]  < \infty \ \Leftrightarrow \exists C > 0, \ \forall x \ge 0 \ \Rightarrow
{\bf P} (|\xi| > x) \le e^{ - C x^{m} }, \ m = \const > 0, \eqno(4.4)
$$
see \cite{Kozachenko1}, \cite{Ostrovsky3}, chapter 1, section 3.\par

\vspace{3mm}

 It follows from proposition (2.3) that if $  B < 2,  $ then the space
$ G(\psi ; (1,B)) $  is not w.n.r.i. space. \par

\vspace{3mm}

{\it  Therefore, we will suppose in what follows that} $  B \ge 2; $
and we will distinguish two cases: $  2 \le B < \infty $   and $  B = \infty. $\par

\vspace{3mm}

 Let us define for arbitrary function $ \psi(\cdot) \in G(\psi; (1,B) ) = G(\psi; (2,B)) $ the new function

$$
\psi_R(p) := R(p) \cdot \psi(p),  \eqno(4.5)
$$
 The symbol "R" in (4.5) appears in the honour of Rosenthal. \par

\vspace{3mm}

{\bf Theorem 4.1.} Let $ Y_{\psi} = G(\psi; (2,B)),  $ where $ B > 2; $  may be $ B = \infty. $  Then  $ Y_{\psi}  $
is w.n.r.i. space with

$$
CLT(Y_{\psi}) = G(\psi_R),
$$
and

$$
\forall \eta: {\bf E} \eta = 0, \ \eta \in Y_{\psi} \ \Rightarrow ||\eta||G(\psi_R) \le ||\eta||G(\psi). \eqno(4.6)
$$

\vspace{3mm}

{\bf Proof.} Suppose the mean zero r.v. $ \eta $ belongs to the space $  G(\psi); $ we can and will  suppose also without loss of generality
$  ||\eta||G(\psi) = 1; $  then $ |\eta|_p \le \psi(p), p \in (2,B). $\par

We deduce using Rosenthal's inequality,  taking into account the restriction $ p \ge 2: $

$$
| n^{-1/2} \sum_{i=1}^n \eta(i)  |_p \le R(p) \cdot |\eta|_p \le R(p) \cdot \psi(p)  = \psi_R(p),
$$
or equally

$$
|| n^{-1/2} \sum_{i=1}^n \eta(i)  ||G(\psi_R) \le 1 = ||\eta||G(\psi),
$$
Q.E.D.\par

\vspace{3mm}

 Let us define

$$
i_{(\psi)}(\theta) = || l(\eta, \theta)||G(\psi) = || \partial \ln g(\eta,\theta)/\partial \theta ||G(\psi),
$$
 the Fisher's information relative the space $  G(\psi). $ We conclude on the basis on the conditions os theorem
 4.1 in the case of a sample of the volume $ n $  for any regular non - biased estimate $ \hat{\theta}_n:  $

$$
\sqrt{n} ||\hat{\theta}_n  - \theta_0||G'(\psi_R) \ge \frac{1}{i_{(\psi)}(\theta_0)}. \eqno(4.7)
$$

 Note that the associate space to the GLS are investigated in the articles \cite{Fiorenza1}, \cite{Fiorenza2},
\cite{Kokilashvili1}, \cite{Ostrovsky4}, \cite{Liflyand1}. \par

 Let us consider the two cases: $  B < \infty $ and $ B = \infty. $\par

\vspace{3mm}

{\bf First case: $ B < \infty. $}\par

\vspace{3mm}

 In this case holds true the simple estimate:

$$
\psi_R(p) \le K(B) \cdot \psi(p) := C \cdot \frac{B}{e \ln B} \cdot \psi(p), \ C = 1.7768\ldots,
$$
so that the norms $ ||\cdot||G(\psi)  $ and $ ||\cdot||G(\psi_R) $ and correspondingly the norms $ ||\cdot||G'(\psi)  $
and $ ||\cdot||G'(\psi_R) $  are equivalent. The inequality (4.7) may be transformed in the considered case as follows:

$$
\sqrt{n} ||\hat{\theta}_n  - \theta_0||G'(\psi) \ge \frac{1}{K(B) \ i_{(\psi)}(\theta_0)}. \eqno(4.8)
$$

 For instance, the function $  \psi(p) $ may be as follows:

$$
\psi(p) \asymp (B - p)^{-\gamma}  \ L(B - p), \  B = \const \ge 2, \  p \in (2,B),
$$
 $  L(x) $ is positive continuous slowly varying as $ x \to 0+ $ function,  see \cite{Ostrovsky4}, \cite{Liflyand1}.\par

 \vspace{3mm}

{\bf Second case: $ B = \infty. $} \par

\vspace{3mm}

 Define the function

$$
\psi_m(p) = p^{1/m}, \ m = \const > 2, \ p \ge 2,
$$
and introduce the following Grand Lebesgue Space $  G_m: $

$$
G_m = G(\psi_m) = \{ \eta: \ {\bf E} \eta = 0, \ ||\eta||_m := \sup_{ p \ge 2} |\eta|_p/p^{1/m} < \infty \}. \eqno(4.9)
$$
 As we know, the {\it centered} r.v. $ \eta $ belongs to this space iff

$$
{\bf P} (|\eta| > x) \le \exp \left( -(x/C(m))^m \right), C(m) = \const > 0.
$$

 Note that the case $ m=2 $ correspondent to the so - called subgaussian variables, see \cite{Buldygin1}, \cite{Buldygin2},
 \cite{Kahane1}, \cite{Kozachenko1}, \cite{Ostrovsky6}.\par
 Let us consider the symmetrical distributed r.v. $ \eta $ such that

$$
{\bf P} (|\eta| > x) = \exp \left( -x^m \right), x > 0;
$$
then $ \eta \in G_m $ and $  0 < C_1(m) \le ||\eta||G_m \le C_2(m) < \infty; $ but evidently

$$
\sup_n || n^{-1/2} \sum_{i=1}^n \eta(i) ||G_m  = \infty.
$$
 This example imply that the space $ CLT(Y) $ may by essentially different from the source r.i. space $ Y. $ \par

 It is easy to verify that in this case the $ CLT(G_m) $ space coincides  with the subgaussian space $  G_2. $ \par

\vspace{4mm}

\section{ Exponential Orlicz's norm }

\vspace{3mm}

 Let $ \phi = \phi(\lambda), \lambda \in (-\lambda_0, \lambda_0), \ \lambda_0 = \const \in (0, \infty] $
be some even strong convex which takes positive values for positive arguments twice continuous
differentiable function, such that
$$
 \phi(0) = 0, \ \phi^{//}(0) > 0, \ \lim_{\lambda \to \lambda_0} \phi(\lambda)/\lambda = \infty. \eqno(5.1)
$$
 We denote the set of all these function as $ \Phi; \ \Phi = \{ \phi(\cdot) \}. $ \par

 We say that the {\it centered} random variable (r.v) $ \xi = \xi(\omega) $
belongs to the space $ B(\phi), $ if there exists some non-negative constant
$ \tau \ge 0 $ such that

$$
\forall \lambda \in (-\lambda_0, \lambda_0) \ \Rightarrow
{\bf E} \exp(\lambda \xi) \le \exp[ \phi(\lambda \ \tau) ]. \eqno(5.2).
$$
 The minimal value $ \tau $ satisfying (4) is called a $ B(\phi) \ $ norm
of the variable $ \xi, $ write
 $$
 ||\xi||B(\phi) = \inf \{ \tau, \ \tau > 0: \ \forall \lambda \ \Rightarrow
 {\bf E}\exp(\lambda \xi) \le \exp(\phi(\lambda \ \tau)) \}. \eqno(5.3)
 $$
 This spaces are very convenient for the investigation of the r.v. having an
exponential decreasing tail of distribution, for instance, for investigation of the limit theorem,
the exponential bounds of distribution for sums of random variables,
non-asymptotical properties, problem of continuous of random fields,
study of Central Limit Theorem in the Banach space etc.\par

  The space $ B(\phi) $ with respect to the norm $ || \cdot ||B(\phi) $ and
ordinary operations is a Banach space which is isomorphic to the subspace
consisting  on all the centered variables of {\it exponential} Orlicz’s space
$ (\Omega,B,{\bf P}), N(\cdot) $ with $ N \ - $ function

$$
N(u) = \exp(\phi^*(u)) - 1, \ \phi^*(u) = \sup_{\lambda} (\lambda u - \phi(\lambda)). \eqno(5.4)
$$
 The transform $ \phi \to \phi^* $ is called Young-Fenchel transform. The proof
of considered assertion used the properties of saddle-point method and theorem of Fenchel-Moraux:
$$
\phi^{**} = \phi.
$$
 The detail investigation of these spaces see in \cite{Kozachenko1}, \cite{Ostrovsky3}, chapters 1,2.
For example, this spaces are a particular cases of $ G(\psi) $  spaces. Namely, if $ \lambda_0 = \infty, $
then the $  B(\phi) $  space is isomorphic to the $ G(\psi_{\phi}) $ space with

$$
\psi_{\phi}(p) = \frac{p}{\phi^{-1}(p)}, \ p \ge 2,
$$
see \cite{Kozachenko1}. \par

 There is a proof also in particular that if the mean zero non - trivial r.v. $ \xi $ belongs to the space  $  B(\phi), $ then

 $$
 \max ( {\bf P}(\xi > x), {\bf P} (\xi < -x)) \le \exp \left(-\phi^*(x/\tau) \right), \ x \ge 0,
 $$
exponential tail estimate; and is true the inverse inequality: if the centered r.v. $  \xi  $ satisfies the last inequality, then
it belongs to the space $ B(\phi): \ ||\xi||B(\phi) \le C(\phi) \cdot \tau. $ \par

 But we can strengthen in  the case of $ B(\phi)  $  spaces some assertions on the  Grand Lebesgue Spaces.  \par

 As in the last section, the function $  \phi(\lambda) $ may be introduced constructively. Namely, let $ \{ \xi(\alpha) \}, \
\alpha \in A  $ be a family of centered r.v., satisfying the uniform Kramer's condition:

$$
\exists C = \const > 0 \ \Rightarrow \sup_{\alpha \in A} {\bf P} (|\xi(\alpha)| > x) \le \exp( - C \cdot x), \ x \ge 0.
$$

 Then we can define

$$
\phi^{(A)}(\lambda) := \sup_{\alpha \in A} \ln {\bf E} \exp (\lambda \xi(\alpha) ), \ |\lambda| < \lambda_0 = \const > 0.
$$

\vspace{3mm}

The associate (and dual) space to the Orlicz  spaces are described, for example, in the famous book of M.M.Rao and Z.D.Ren
\cite{Rao1}, chapter 4,5.\par

 Let $ \phi(\cdot) \in \Phi; $ define a new function

$$
\overline{\phi}(\lambda) = \sup_{n = 1,2,\ldots} [n \cdot \phi(\lambda/\sqrt{n})]; \eqno(5.5)
$$
 then $ \overline{\phi}(\cdot) \in \Phi. $ \par

For example, if $ \phi(\lambda) \asymp \lambda^Q,  \ \lambda \ge 1, \ Q = \const > 1, $  then
$ \overline{\phi}(\lambda) \asymp \lambda^{\max(2,Q)}, \ \lambda \ge 1. $ \par

 Notice that if $ Q \ge 2, $ then $ \overline{\phi}(\lambda) \asymp \phi(\lambda), \ \lambda \ge 1. $ This possibility
is absent for the GLS spaces.\par

 The subgaussian r.v. forms the $ B(\phi_2) = B_2 $  space  with $ \phi_2(\lambda) = 0.5 \lambda^2, \ \lambda \in R;  $ in this case
$ \overline{\phi_2}(\lambda) = \phi_2 (\lambda).  $ \par

\vspace{3mm}

{\bf Theorem 5.1.} Let $ Y_{\phi} = B(\phi),  $ where $ \phi \in \Phi. $   Then  $ Y_{\phi}  $ is w.n.r.i. space with

$$
CLT(Y_{\phi}) = B(\overline{\phi})
$$

and

$$
\forall \eta: {\bf E} \eta = 0, \ \eta \in Y_{\phi} \ \Rightarrow ||\eta||B(\overline{\phi}) \le ||\eta||G(\psi). \eqno(5.6)
$$

\vspace{3mm}

{\bf Proof} is complete analogously to one in theorem 4.1; it based on the equality

$$
{\bf E} \exp \left(\lambda n^{-1/2} \sum_{i=1}^n \eta(i)   \right) \le \exp \left( \overline{\phi}(\lambda)  \right)
$$
or equally

$$
|| n^{-1/2} \sum_{i=1}^n \eta(i) ||B(\overline{\phi}) \le ||\eta||B(\phi). \eqno(5.7)
$$

\vspace{3mm}

 Let us define as before

$$
i_{(\phi)}(\theta) = || l(\eta, \theta)||B(\phi) = || \partial \ln g(\eta,\theta)/\partial \theta ||G(\phi),
$$
 the Fisher's information relative the space $  B(\phi). $ We conclude  under  the conditions os theorem 5.1 in the case of a
sample of the volume $ n $  for any regular non - biased estimate $ \hat{\theta}_n:  $

$$
\sqrt{n} ||\hat{\theta}_n  - \theta_0||B'(\overline{\phi}) \ge \frac{1}{i_{(\phi)}(\theta_0)}. \eqno(5.8)
$$

\vspace{4mm}

\section{ Lorentz norm }

\vspace{3mm}

  Recall that the norm of a r.v. $  \zeta $ in the Lorentz space $ L_{p,q} = L_{p,q}(\Omega),  $ more exactly, quasinorm
$  ||\zeta||^*_{p,q}, \ 1 \le p,q \le \infty $   is defined  as follows:

$$
||\zeta||^*_{p,q} \stackrel{def}{=} \left( \int_0^{\infty} \left[ {\bf P} \{ |\zeta| \ge x  \} \right]^{q/p} \ d x^q  \right)^{1/q}, \ 1 \le p,q < \infty,
$$
and

$$
||\zeta||^*_{p,\infty} \stackrel{def}{=} \sup_{x > 0} \left[ x \left( {\bf P} \{ |\zeta| \ge x  \} \right)^{1/p} \right].
$$

 The detail investigation of these spaces see in the books \cite{Bennet1}, \cite{Stein1}; for instance, it is proved that this quasinorm is linear
equivalent to really norm

$$
||\zeta||_{p,q} := \sup_{A: {\bf P}(A) > 0} \left[ \frac{\int_A |\zeta(\omega)| \ {\bf P}(d \omega)}{\nu_{p,q}( {\bf P}(A) )} \right]
$$
for some positive for positive values $ z  $ function $ \nu_{p,q}(z), \ z \in (0,1); \ \nu_{p,q}(0) = 0. $ \par
 Using for us important facts about these spaces are obtained the book of M.Sh.Braverman \cite{Braverman1}.\par
  In particular, $ L_{p,p} = L_p,  $ therefore the Lorentz are direct generalization  of Lebesgue - Riesz spaces. But the exact values
of Rosenthal's constants for this spaces are now unknown.\par

\vspace{3mm}

{\bf  Theorem 6.1.}   Denote $ r = \min(p,q). $  The Lorentz space $  L_{p,q}(\Omega) $ is strong
normal r.i. space iff $  r > 2 $ or $  q = 2 \le p. $\par

\vspace{3mm}

 This assertion is in fact proved in the book of M.Sh.Braverman  \cite{Braverman1}, p. 11 - 13, theorem 7.\par

\vspace{4mm}

 \section{ Concluding remarks }

 \vspace{3mm}

 The function $ \psi = \psi(p) $ and $ \phi = \phi(\lambda)  $ in the sections 4 and 5 may be constructively introduced.  Indeed,
in the case of $ G(\psi) $ spaces we can introduce the {\it natural} function for the  family of the r.v. $ l(\eta, \theta), \ \theta \in \Theta:  $

$$
\psi_0(p) := \sup_{\theta \in \Theta} \left[ \ \int_X  |g'_{\theta}(x,\theta)|^p \ g^{1 -p}(x,\theta) \ \mu(dx) \  \right]^{1/p}, \eqno(7.1)
$$
if of course the last expression is finite  for some values $ p $ greatest than 2.\par

 Notice that this choice of this function $ \psi_0(p) $ is {\it optimal,} i.e. minimal.\par

 If for instance the density $ g(x,\theta) $ has a form $  g(x,\theta) = g_0(x-\theta), \ x,\theta \in R $ ("shift" case),
 where $ g_0(\cdot) $ is differentiable density function,
 and $ \mu(A) $ is an ordinary  Lebesgue measure, then the integral inside the expression (7.1) does not depended on the value
$ \theta $ and hence

$$
\psi_0(p) := \left[ \ \int_R  |g'_0(x)|^p \ g_0^{1 - p}(x) \  dx \  \right]^{1/p}, \eqno(7.2)
$$
if of course $ \exists p_0 > 2, \ \psi_0(p_0) < \infty. $ \par

\vspace{3mm}

 Another example - scaling parameter. Here again $ X = R, \ \mu $ is Lebesgue measure, but $ \theta \in (\theta_-, \theta_+), \ \theta_- > 0 $
 is scaling parameter:
 $$
 g(x,\theta) =  \theta^{-1} h(x/\theta),
 $$
where $ h(\cdot) $ is differentiable density function. Then

$$
 \left[ \ \int_X  |g'_{\theta}(x,\theta)|^p \ g^{1 -p}(x,\theta) \ \mu(dx) \  \right]^{1/p} =
 \theta^{-1} \cdot \left[ \int_R |h(y) + y h'(y)|^p \cdot h^{1 - p}(y) \ dy \right]^{1/p}. \eqno(7.3)
$$

\vspace{3mm}

The correspondent natural function $ \phi_0(\lambda) $  for the  family  $  \{ l(\eta, \theta) \}  $ in the space $  B(\phi) $
look not so nice as for the Grand Lebesgue Spaces:

$$
\phi_0(\lambda) := \sup_{\theta \in \Theta} \log \int_X \exp \left( \lambda \frac{g'_{\theta}(x, \theta) }{g(x,\theta)} \right)
 \ g(x, \theta) \ \mu(dx), \eqno(7.4)
$$
if it is finite  for sone non - trivial interval $ |\lambda| < \lambda_0, \ \lambda_0 = \const > 0; $ with evident modification
for the shift or scale parameter $  \theta. $\par

\vspace{3mm}

\begin{center}

{\bf  Conclusions: }\par

\end{center}

\vspace{3mm}

 The assertions of the sections 2-6 may be simplify as follows: under appropriate conditions on the {\it weak} r.i. space $ (Y, ||\cdot||Y),  $
for example the space $ L_q(\Omega,{\bf P}) $ with $ 1 < q < 2, $
 and on the density $ g(\cdot,\cdot), $ for {\it arbitrary} unbiased regular sample estimate $ \hat{\theta}_n $ holds true the following inequality:

$$
\sqrt{n} \cdot || \hat{\theta}_n  - \theta_0||Y  \ge C(Y, g(\cdot),\theta_0 ), \ n = 1,2,\ldots; \eqno(7.4)
$$
while for the {\it MLE} sample estimate  $ \tilde{\theta}_n $ and for {\it  strong} r.i. space $ (Z, ||\cdot||Z),  $
for example, for $  L_p(\Omega,{\bf P}) $ with $  p > 2,  \ G(\psi) \ $ and $ B(\phi) $ spaces is true the opposite inequality:

$$
\sqrt{n} \cdot || \tilde{\theta}_n  - \theta_0||Z  \le \tilde{C}(Z, g(\cdot, \cdot), \theta_0 ), \ n = 1,2,\ldots. \eqno(7.5)
$$

\end{document}